\documentclass[singlespacing]{bmcart}
\usepackage[T1]{fontenc}
\usepackage[pdftex]{graphicx}

\usepackage[cmex10]{amsmath}

\usepackage{eurosym}

\usepackage{hyperref}
\usepackage{booktabs}

\usepackage{tikz}
\usepackage{pgfplots}
\pgfplotsset{compat=1.17} 
\usetikzlibrary{positioning,shapes.geometric}
\pgfdeclarelayer{background}
\pgfdeclarelayer{foreground}
\pgfsetlayers{background,main,foreground}

\usepackage[noend]{algorithm2e}
\SetKw{And}{and}{}{}

\usepackage{mathtools}
\usepackage{amssymb}
\usepackage{subfig}
\usepackage{color, colortbl}
\newcommand{\Vfinal}{{V_{\mathrm{final}}}}
\newcommand{\Vinit}{{V_{\mathrm{init}}}}
\newcommand{\R}{\mathbb{R}}
\newcommand{\N}{\mathbb{N}}

\newcommand{\Mod}[1]{\ \mathrm{mod}\ #1}
\newcommand{\etain}{\eta_\mathrm{in}}
\newcommand{\etaout}{\eta_\mathrm{out}}
\DeclareMathOperator*{\argmin}{arg\,min}

\newcommand{\new}[1]{\textcolor{black}{#1}}

\usepackage[utf8]{inputenc} 

\startlocaldefs
\endlocaldefs

\begin{document}

\begin{frontmatter}

\begin{fmbox}
\dochead{Research}


\title{Determining Cost-Efficient Controls of Electrical Energy Storages Using Dynamic Programming}


\author[
  addressref={aff1},                   
  email={stiglmayr@uni-wuppertal.de}   
]{\inits{M.S.}\fnm{Michael} \snm{Stiglmayr}}
\author[
  addressref={aff1},
  email={suhlemeyer@uni-wuppertal.de}
]{\inits{S.U.}\fnm{Svenja} \snm{Uhlemeyer}}
\author[
  addressref={aff2},                   
  email={bjoern.uhlemeyer@uni-wuppertal.de}   
]{\inits{B.U.}\fnm{Björn} \snm{Uhlemeyer}}
\author[
  addressref={aff2},
  email={zdrallek@uni-wuppertal.de}
]{\inits{M.Z.}\fnm{Markus} \snm{Zdrallek}}


\address[id=aff1]{
  \orgdiv{School of Mathematics and Natural Sciences, IZMD},             
  \orgname{University of Wuppertal},          
  \city{Wuppertal},                              
  \cny{Germany}                                    
}
\address[id=aff2]{%
  \orgdiv{School of Electrical, Information and Media Engineering},
  \orgname{University of Wuppertal},
  \city{Wuppertal},
  \cny{Germany}
}





\begin{abstractbox}

\begin{abstract} 
Volatile electrical energy prices are a challenge and an opportunity for small and medium-sized companies in energy-intensive industries. 
By using electrical energy storage and/or an adaptation of production processes, companies can significantly profit from time-depending energy prices and reduce their energy costs.

We consider a time-discrete optimal control problem to reach a desired final state of the energy storage at a certain time step. Thereby, the energy input is discrete since only multiples of 100 kWh can be purchased at the EPEX SPOT market. 
We use available price estimates to minimize the total energy cost by a rounding based dynamic programming approach.
With our model non-linear energy loss functions of the storage can be considered and we obtain a significant speed-up compared to the integer (linear) programming formulation.
\end{abstract}




\end{abstractbox}
\end{fmbox}

\end{frontmatter}



\allowdisplaybreaks

\section{Introduction}\label{sec:intro}
The climate targets for Germany to become climate-neutral by 2045, which were reinforced in 2021, imply that the ambitions for expanding renewable energy must be further increased. As a result, the share of volatile power generation will rise and energy storage will become increasingly important. In 2020, there were already almost 300 hours with negative electricity prices on the day-ahead market~[\url{www.epexspot.com}]. This results, among other things, from the oversupply by renewable power generation. The reversal of the trend for an increasing number of negative electricity prices can be significantly influenced by energy storage. Energy storage systems can assume different functions in the energy system. In multi-use approaches, it is increasingly being investigated how battery storage can be used in times when it is not required for its primary purpose, such as primary control power generation or self-supply optimization \cite{Zeh.2016}.
One option is basically to trade energy on EPEX Spot and take advantage of the price spread between two points in time. Besides mixed-integer optimization models \cite{Kumtepeli.2020}, various (approximate) dynamic programming approaches (see, e.g., \cite{loehndorf10optimal,salas18bench,goubko22approximate,scott22approximate,ruether22iterative}) are used to determine cost efficient controls of electrical storages and/or grids. 
Most approaches therefore use \textit{approximate dynamic programming} \cite{powell11approximate} in order to avoid the ``curses of dimensionality'' by approximating the value function in each state. In \cite{halman18complexity} it is shown that the energy storage problem can be solved in polynomial time in a deterministic setting, while it is an \textsf{\textbf{NP}}-complete problem if pri\new{c}es and energy production are stochastic.

In this paper we consider a time-discrete optimal control problem of an electrical energy storage device and present a rounding based dynamic programming approach, which considers a discrete state space by rounding the energy level in the storage. This  also reduces the computation time significantly in contrast to the solution of the mixed-integer programming problem. 
The rounding of the state space enables us to optimize the control over longer time periods (up to one year), which is of particular interest for the layout of the storage device within a retrospective analysis.

The calculations shown here are based on the example of electro-chemical energy storage. In principle, however, these considerations apply to all forms of energy storage (electrical, electromagnetic, electro-chemical, mechanical, thermal and chemical) and the suggested algorithm can be applied analogously.

This paper is organized as follows. In Section~\ref{sec:prob} we present a time-discrete model of an electrical energy storage devise which takes energy loss during charging, withdrawal, and self-discharge into account. Based on this model we develop a mixed-integer programming (MIP) formulation for the optimal control of the energy storage. Due the limited applicability of the MIP for larger instances we present a rounding based dynamic programming approach in Section~\ref{sec:dyn}, which efficiently approximates the problem. Numerical tests are presented in Section~\ref{sec:num}, including a run-time comparison of the considered models and a trade-off analysis of the investments in storage devices. Section~\ref{sec:concl} concludes the article and gives a brief outlook on future research directions based on our approach.

\section{Modeling of an Electrical Energy Storage}\label{sec:prob}
\begin{figure}[t]
    \centering
    \captionsetup{type=figure}
    \includegraphics[width=0.7\textwidth]{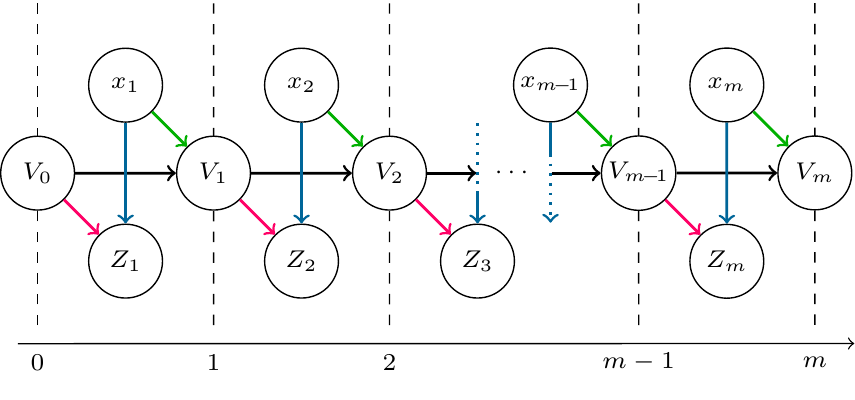}
    \caption{Schematic visualization of an electrical energy storage with fill level $V$ for discretized time intervals. The energy flow is depicted by arrows of different colors. The left energy from previous time intervals (\emph{black}) remains in the storage, possibly affected by some energy loss. The purchased energy $x$ can either be fed into the storage (\emph{green}) or directly consumed (\emph{blue}). The required amount of energy $Z$ can also be (partly) extracted from the electrical storage (\emph{red}).} 
    \label{fig:graph}
\end{figure}

The energy market has varying energy prices due to supply and demand reasons. On account of, \emph{e.g.}, solar or wind energy, prices are relative to the weather, \emph{i.e.,} low costs correlate with using more eco-power. Therefore, storing energy instead of always buying exactly the required amount may be economically as well as ecologically reasonable (see, \textit{e.g.}, \cite{paulus11potential,meese15opt,meese17intra})\new{.}
A comprehensive introduction to dynamical energy pri\new{c}es and their impact on industrial processes is given in \cite{meese18dynamische}.

\subsection{Linear Programming}\label{subsec:lp}
Regarding the day-ahead market, trading is only possible at fixed points in time, hence we define the set of discrete trading dates as $T\coloneqq\{1,\ldots,m\}$, $m\in\N$. Let then $V_t$ denote the charge level of an electrical storage at the time of $t\in T$, which is restricted by lower and upper capacity bounds $c,C\in\R_{\geq0}$, \textit{i.e.}, \(c\leq V_t\leq C\) for all \(t\in T\). Furthermore, let the initial and a final fill level be given, denoted by $V_0 = \Vinit$ and $V_m \geq \Vfinal$, respectively (see Figure~\ref{fig:graph} for a schematic illustration). The lower and upper bounds on the purchased energy per time step are denoted by $l$ and $u$ (\(0\leq l \leq u\)), \emph{i.e.,} $l \leq x_t \leq u$ for all $t\in T$. The fill level $V_t$ depends on three quantities, firstly on the purchased energy $x_t$, secondly on the external energy consumption $Z_t\in\R_{\geq0}$ and lastly on the previous fill level $V_{t-1}$. Regarding the latter, we introduce an energy loss function ${g:\R_{\geq0}\to\R_{\geq0}}$, which is often assumed to be linear, \emph{e.g.}, $g(V) = (1-\beta)\, V$ for some given value $\beta\in(0,1)$. 

The efficiency of storing to and withdrawing energy from the storage is modeled by the efficiency factors $\etain, \etaout\in[0,1]$, which are either assumed to be constant factors, or depending on the amount of stored/withdrawn energy $\etain, \etaout: \R\to[0,1]$. Consequently, it has to be distinguished whether the consumed energy $Z_t$ is taken from the storage or purchased energy is used directly. Let $y_t\in[0,x_t]$ be the amount of energy stored in time step \(t\) and $\zeta_t \coloneqq Z_t - x_t + y_t$ the amount of energy loaded from the storage. In total, we obtain the fill level $V_t$ by
\new{
\begin{equation*}
   V_t = \etain \cdot y_t + g(V_{t-1}) - \frac{1}{ \etaout}\,\zeta_t .
\end{equation*}
To illustrate this formula we consider the two extreme cases: If the energy consumption \(Z_t\) in time step \(t\) is taken completely from the storage (since there is not energy input in this time step, i.e., \(x_t=y_t=0\)) then the energy level in the storage \(V_t=g(V_{t-1}) - \frac{1}{ \etaout}\,Z_t\) is reduced by \(\frac{1}{ \etaout}\,Z_t\) taking into account the energy loss when withdrawing energy from the storage. On the other hand, if the energy consumption equals the energy input \(Z_t=x_t\) energy level of the storage is \(V_t=g(V_{t-1})\) unchanged apart from the time-dependent energy loss \(g\).}

Forecasting models can provide a prognosis about the energy prices $p_t$, $t\in T$ for the period $T$. For a retrospective analysis, however, we can also consider the true prices, \textit{e.g.}, to evaluate the capacity of the energy storage. In the following we will concentrate on the running energy costs and neglect acquisition and other types of fixed costs. Thus, we aim to minimize the total costs of purchased energy $x = (x_1,\ldots,x_m)$, \emph{i.e.,}
    $\min\;\, \sum_{t\in T} p_t \, x_t$. 
Together with the aforementioned constraints, we formulate the following linear program (LP):
\begin{subequations}  
\begin{alignat}{3}\label{eq:lp}
    \min\quad  &\sum_{t\in T}  p_t \, x_t  \\
    \text{s.t.} \quad &l \leq  x_t  \leq u &&\forall t\in T\label{eq:first_nb}\\
    &c \leq  V_t  \leq C &&\forall t\in T\\
    &V_0 = \Vinit\\
    &V_t  = \etain\cdot y_t + g(V_{t-1}) - \frac{1}{\etaout}\, \zeta_t \quad && \forall t\in T\\
    &0 \leq y_t \leq x_t && \forall t\in T\\
    &Z_t = (x_t-y_t) +  \zeta_t && \forall t\in T \\
    &V_m  \geq \Vfinal \\  
    & \new{\zeta_t \geq 0} && \new{\forall t \in T}\\
    & \new{V_t \geq 0} && \new{\forall t \in T}\\
    &\new{(y_1, \ldots, y_m)  \in \R^m \label{eq:last_nb}}\\
    &\new{(x_1,\ldots,x_m) \in \R^m\label{eq:variables} }
\end{alignat}
\end{subequations}
It is well known that linear optimization problems are efficiently solvable (in polynomial time with interior point methods, cf.~\cite{gondzio12interior}). However, the energy is often traded in discrete quantities, which makes the energy input \(x\) a discrete variable. For example, at the EPEX SPOT market only multiples of $100$\,kWh can be purchased. We thus obtain a mixed-integer (linear) programming problem (MIP).

\subsection{Integer Programming}\label{subsec:ip}
Since the energy input can only attain discrete values, we modify equation (\ref{eq:variables}) in (MIP) \new{to equation (\ref{eq:discrete_variables})}: 
\begin{subequations}\label{eq:ip}
\begin{align}
    \min\quad  &\sum_{t\in T}  p_t \, x_t  \\
    \text{s.t.} \quad &\eqref{eq:first_nb}-\eqref{eq:last_nb}\\
    &(x_1,\ldots,x_m) \in h_x\cdot\N^m \; . \label{eq:discrete_variables}
\end{align}
\end{subequations}
Thereby, $h_x\in\R_{\geq 0}$ denotes the discretization step size of the energy input, \emph{i.e.,} \(x\) is restricted to multiples of $h_x$: \(x\in\{0,h_x,2\,h_x,\ldots\}\). 
In contrast to LPs, MIPs 
are in general \textsf{\textbf{NP}}-hard optimization problems, which are solvable by a significant computational effort, \textit{e.g.}, using branch and bound based approaches \cite{wolsey98integer}.
Moreover, if the energy loss in the storage depends non-linearly on the fill level, one would obtain a mixed-integer non-linear optimization problem, which are computationally even more demanding.

On that account, we introduce our dynamic programming approach in the following chapter.

\section{Rounding-based Dynamic Programming}\label{sec:dyn}
\subsection{Dynamic Programming}

The central idea of dynamic programming is to break down an optimization problem into a sequence of smaller efficiently solvable subproblems. Thereby, dynamic programming relies on \emph{Bellman's principle of optimality} \cite{bellman57dynamic}, which states that a solution can only be optimal if its intermediate solutions (up to a certain state/time) are optimal w.r.t.\ the corresponding subproblems. 
Knapsack problems \cite{kellerer04knapsack} and shortest path problems are the most prominent examples of optimization problems satisfying Bellman's principle \cite{papadimitriou82combi}, which does not hold for all optimization problems. 
The discrete electrical energy storage problem \eqref{eq:ip}, however, satisfies Bellman's principle, since a control \(x_t,y_t\), \(t=1,\ldots,\tau\), of the storage up to an intermediate time step \(\tau\in T\) with fill level of \(V_\tau\) can only be extended to an optimal solution if there is no other feasible policy reaching this (or a larger) fill level \(V_\tau\) at time \(\tau\) with lower energy cost. 
The optimal solution of the overall problem can then be derived from the optimal solutions of these subproblems. 
Applying Bellman's recursion \cite{bellman57dynamic} we determine the cheapest way to reach every feasible fill level at time step \(t\) based on the costs at time step \(t-1\). The optimal control $(x^\ast_1,\ldots,x^\ast_m)$ for an arbitrary final fill level can then be reconstructed by a backtracking procedure.
Adapted to the previously introduced electrical storage problem \eqref{eq:ip}, we initialize the recursion for the \new{total energy costs \(z_t(d)\) up to time step \(t\) to reach a given storage fill level \(d\)} as
\begin{equation*}
    z_1(d) \coloneqq 
    \begin{cases}
        p_1 \cdot x_1(d)  & \text{, if } x_1(d) \Mod h_x = 0\\
        \infty & \text{, otherwise}    
    \end{cases} \; ,
\end{equation*}
where $x_1(d) = (d - g(\Vinit) + \new{\frac{\zeta_1}{\etaout}})/\etain $ is the amount of energy required to reach the level of \(d\) in the current state. \new{Thereby, we assume that energy is only withdrawn from the storage when it is not necessary to reach the desired storage level \(d\) in time step \(t\), i.e., \(\zeta_t=\max\{0,\etaout(g(V_{t-1})-d)\}\), since it is always preferable to directly consume energy over its lossy storage. Consequently, in each time step there can be only either charging of the storage or withdrawel of energy from the storage.}

\subsection{Rounding in the State Space}
Since the computational efficiency of dynamic programming algorithms strongly depends on the size of the state space a straightforward application of the Bellman recursion onto the discrete storage problem would lead to numerical difficulties. Due to the energy loss function \(g\) it is very unlikely that different policies end up at the same fill level, so the number of states grows exponentially with increasing number of time steps.
In order to limit the number of states, we discretize the state space, \emph{i.e.,} the fill level of the storage, with step size $h_V\in\R_{\geq 0}$. A similar approach is proposed in \cite{luus90optimal} on continuous control problems.

Let
\begin{equation*}
    \lceil V \rceil_{h} \coloneqq \begin{cases} V & \mathrm{, if }~ V \Mod h = 0 \\ V - (V \Mod h) + h & \mathrm{, otherwise}\end{cases}
\end{equation*}
define the ceil function with respect to some step size $h$,
\begin{equation*}
    \lfloor V \rfloor_{h} \coloneqq V - (V \Mod h)
\end{equation*}
the corresponding floor function, $V\in\R_{\geq 0}$. In the \textit{rounding based dynamic programming (RBDP)} algorithm (Algorithm~\ref{alg:dyn_round2}), we underestimate the fill level of the storage, \emph{i.e.,} we round off in the state space and consider the following adapted recursion formula, by which we obtain the optimal control to approximately reach a given fill level $d$ at time step $t$ based on the controls up to the previous time step \(t-1\): 
\begin{align}
 z_t(d)\coloneqq \min\limits_{d'\in h_V\cdot\N_0} \biggl\{z_{t-1}(d') + p_t\, x_t(d) \colon 
    \Bigl\lfloor \etain\, y_t + g(d') - \frac{\zeta_t}{\etaout} \Bigr \rfloor_{h_V} = d \biggr\}. 
\end{align}
Using this recursion formula we can apply a dynamic programming scheme to compute the optimal control for each energy level of the storage \(d=c,c+h_V,\ldots, C\) for all time steps \(\new{t}=1,\ldots,m\), \new{and each energy input \(k\), by eliminating dominated states, }cf.~Algorithm~\ref{alg:dyn_round2}. \new{In the following we assume that \(l\) and \(u\) are multiples of the discretization \(h_x\).}

\begin{algorithm}[tbh]
  \caption{Rounding based Dynamic Programming \label{alg:dyn_round2}\phantom{\rule{1pt}{4ex}}}
  \KwData{\({h_x,h_V\geq 0, m,\Vinit\new{\in\mathbb{N}, \Vfinal \in h_V\mathbb{N},} c,C,l,u \geq 0}\),
  \({\beta,\new{\etain,\etaout} \in[0,1], Z\in[0,\infty)^m, 
  x:\mathbb{N}\to[l,u]^m}\),
  \({z:\mathbb{N}\to(\mathbb{R}\cup\infty)^m, A:\mathbb{N}\to(\mathbb{N}\cup\infty)^m}\)
  }
  
  \KwResult{\(x^\ast\in[l,u]^m, V\in[c,C]^m\)}

  \For{\(d\coloneqq c,c+h_V,\ldots,C\)}{
    \(A_1(d)\coloneqq \infty\)\;
    \If{\(((d - g(\Vinit) + \zeta_1/\etaout)/\etain) \Mod h_x = 0\) \And \(\new{(d - g(\Vinit) + \zeta_1/\etaout)/\etain <ub }\)}{
      \(x_1(d)\coloneqq (d - g(\Vinit) + \zeta_1/\etaout)/\etain \) \;
      \(z_1(d)\coloneqq x_1(d) \cdot p_1\)\;
    }
    \Else{
      \(z_1(d)\coloneqq \infty\)\;
      \(x_1(d)\coloneqq 0\)\;
    }
  }
  \For{\(t\coloneqq 2,\ldots,m\)}{
    \For{\(d\coloneqq c,c+h_V,\ldots,C\)}{
        \(z_t(d) \coloneqq \infty\)\;
        \(x_t(d)\coloneqq 0\)\;
        \(A_t(d) \coloneqq \infty\)\;
        \For{\(k\coloneqq l,l+h_x,\ldots,u\)}{
          \(\zeta_t \coloneqq \max\{Z_t-k,0\}\)\;
          \(y_t\coloneqq \max\{k-Z_t, 0\}\)\;
          \(lb := \lceil g^{-1}( d + Z_t - \etain \, y_t + \zeta_t/\etaout) \rceil_{h_V}\)\;
          \(ub := \lceil g^{-1}(d + Z_t - \etain \, y_t + \zeta_t/\etaout + h_V) \rceil_{h_V}\)\;
          \For{\(W\new{\coloneqq} lb,\new{lb + h_V,}\ldots,ub\)}{
            \If{\(c \leq W \leq C\) \And \(z_{t-1}(W)+k\cdot p_t < z_t(d) \)}{
                \(z_t(d) \coloneqq z_{t-1}(W)+k\cdot p_t\)\;
                \(x_t(d)\coloneqq k\)\;
                \(A_t(d) \coloneqq W\)\;
            }
          }
        }
    }
}
  \(d^\ast \coloneqq \argmin_{d\in\{\Vfinal,\ldots,C\}} z_m(d)\)\;
  \(z^\ast \coloneqq z_m(d^\ast)\)\;
  \(x^\ast_m \coloneqq x_m(d^\ast)\)\;
  \(V_m \coloneqq d^\ast\)\;
  \For{\(i=m-1,\ldots,1\)}{
    \(d^\ast \coloneqq A_{i+1}(d^\ast)\)\;
    \(x^\ast_i \coloneqq x_i(d^\ast)\)\;
    \(V_i \coloneqq d^\ast\);
  }
\end{algorithm}

By doing so, the rounding error 
\begin{align*}
    \varepsilon \coloneqq \bigl|V - \lfloor V \rfloor_{h_V}\bigr| &= \bigl| V - (V - (V \Mod h_V)) \bigr| 
    = (V \Mod h_V) < h_V
\end{align*}
is bounded from above by $h_V$ for one single time step.
In order to estimate the total error of our algorithm, we first consider some fill level 
\begin{align}
    V_t &= \Bigl\lfloor \etain\, y_t  + g(V_{t-1}) - \frac{1}{\etaout}\,\zeta_t \Bigr\rfloor_{h_V}     
    =\etain\, y_t  + g(V_{t-1}) - \frac{1}{\etaout}\,\zeta_t - \varepsilon_t
\end{align}
with rounding error $0 \leq \varepsilon_t < h_V$, $t\in T$. If we assume that the energy loss function \(g\) is a monotonically increasing, linear function, we can derive an explicit bound on the total rounding error:
\begin{align*} 
    V_m= g^m(\Vinit) + \sum_{i=1}^m g^{m-i}\Bigl(\etain\, y_i - \frac{1}{\etaout}\,\zeta_i - \varepsilon_i\Bigr)\;,
\end{align*}
where \(g^k\) denotes the \(k\)-times iterated function, \textit{i.e.}, \(g^0=\mathrm{id}\) and \(g^k=g(g^{k-1})\) with \(k\in\N\).
Then, the total rounding error is given by
\[
    \varepsilon_\mathrm{tot}\coloneqq\sum_{i=1}^m \Bigl|g^{m-i}\Bigl(\etain\, y_i - \frac{\zeta_i}{\etaout} - \varepsilon_i\Bigr)-g^{m-i}\Bigl(\etain\, y_i - \frac{\zeta_i}{\etaout} \Big)\Bigr|\;.
\]
\new{By definition of the function $g$ as $g(V)=(1-\beta)\, V$ with $\beta\in(0,1)$, it holds that}
\begin{align}
    \varepsilon_\mathrm{tot} &= \sum_{i=1}^m g^{m-i}\Bigl(\etain\, y_i - \frac{\zeta_i}{\etaout} \Bigr) - g^{m-i}\Bigl(\etain\, y_i - \frac{\zeta_i}{\etaout} - \varepsilon_i \Big) \notag\\
    &< \sum_{i=1}^m g^{m-i}\Bigl(\etain\, y_i - \frac{\zeta_i}{\etaout} \Bigr)  - g^{m-i}\Bigl(\etain\, y_i - \frac{\zeta_i}{\etaout} - h_V \Big) \notag\\
    &= \sum_{i=1}^m g^{m-i}(h_V) \leq  m\cdot h_V \;.\label{eq:linear_bound}
\end{align}
Since the rounding based dynamic programming algorithms is an exact method on the discretized state space, the difference between the costs of a solution of RBDP and the cost of the optimal solution of \eqref{eq:ip} are at most \(m\, h_V\cdot \max_t\{p_t\}\).
The bound in equation \eqref{eq:linear_bound} can be computed and subtracted from the upper capacity bound $C$ in order to guarantee feasibility of the exact solution. Hence, $h_V$ should be chosen depending on the maximal capacity, \emph{i.e.,} $\lfloor C - \varepsilon_\mathrm{tot} \rfloor_{h_V} \gg 0$.

\section{Numerical Tests}\label{sec:num}
\begin{table}[t]
    \centering
    \vspace{0.75cm}
    \begin{tabular}{l@{\extracolsep{2ex}}*{4}{@{\extracolsep{4ex}}c}}
        \toprule
        \textbf{Capacity}  & 500 & 1000 & 2500 & 5000  \\\midrule
        \textbf{LP solution} & 13,358\euro\  & 13,224\euro\  & 12,954\euro\  & 12,724\euro\  \\\midrule
        \textbf{MILP solution} & 13,388\euro\  & 13,242\euro\  & 12,980\euro\  & 12,738\euro\ \\\midrule
        \textbf{DP solution} & 13,396\euro\  & 13,250\euro\  & 12,985\euro\  & 12,742\euro\  \\\bottomrule\smallskip
    \end{tabular}
    \caption{Comparison of different optimization algorithms to minimize the energy costs over one week (06/15/18--06/21/18) for four different storage capacities. Purchasing exactly the required amount of energy in each time step, \emph{i.e.,} without using an electrical storage, costs 13,532\,\euro. Note that the LP solutions are only given as lower bounds and do not correspond to feasible controls.}
    \label{tab:comparison}
\end{table}
\begin{figure*}[t]
    \centering
    \captionsetup[subfloat]{labelformat=empty}
    \subfloat[]{{\includegraphics[width=0.95\textwidth]{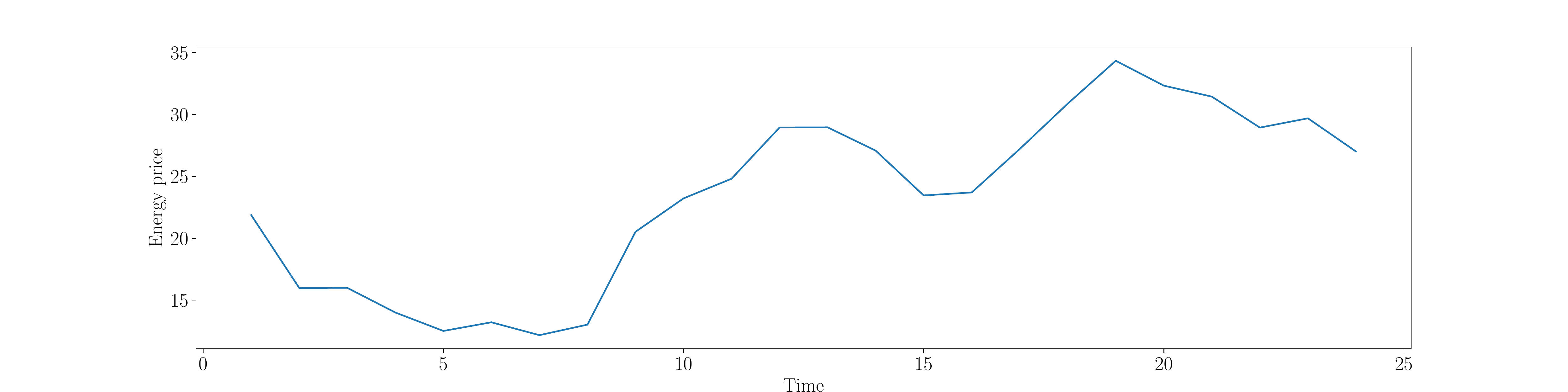}}}\\[-1.5em]
    \captionsetup{type=figure}
    \subfloat[(LP) obj.~value: 1108\,\euro\ ]{\includegraphics[width=0.32\textwidth]{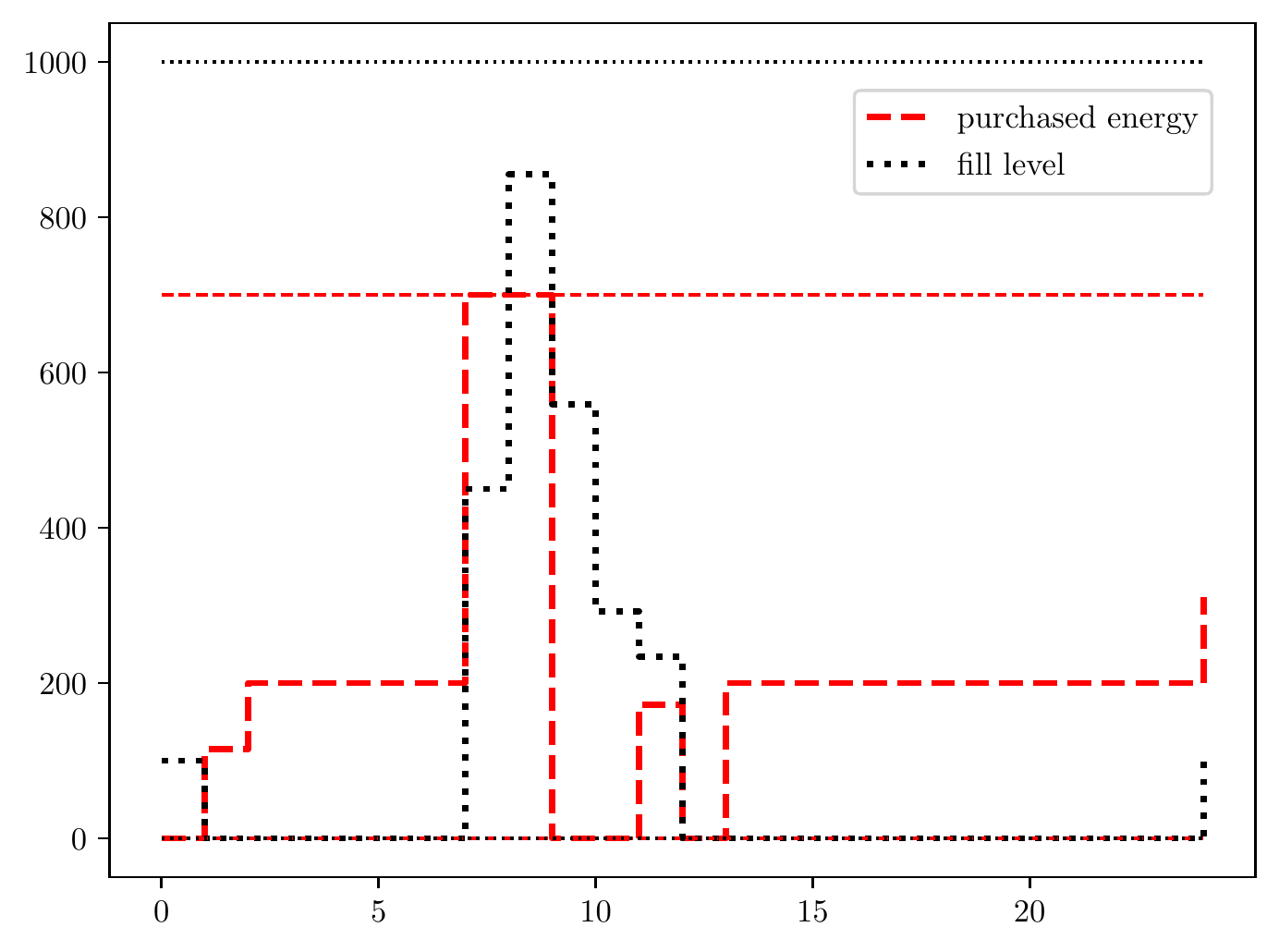}}~
    \subfloat[(MILP) obj.~value: 1126\,\euro\ ]{\includegraphics[width=0.32\textwidth]{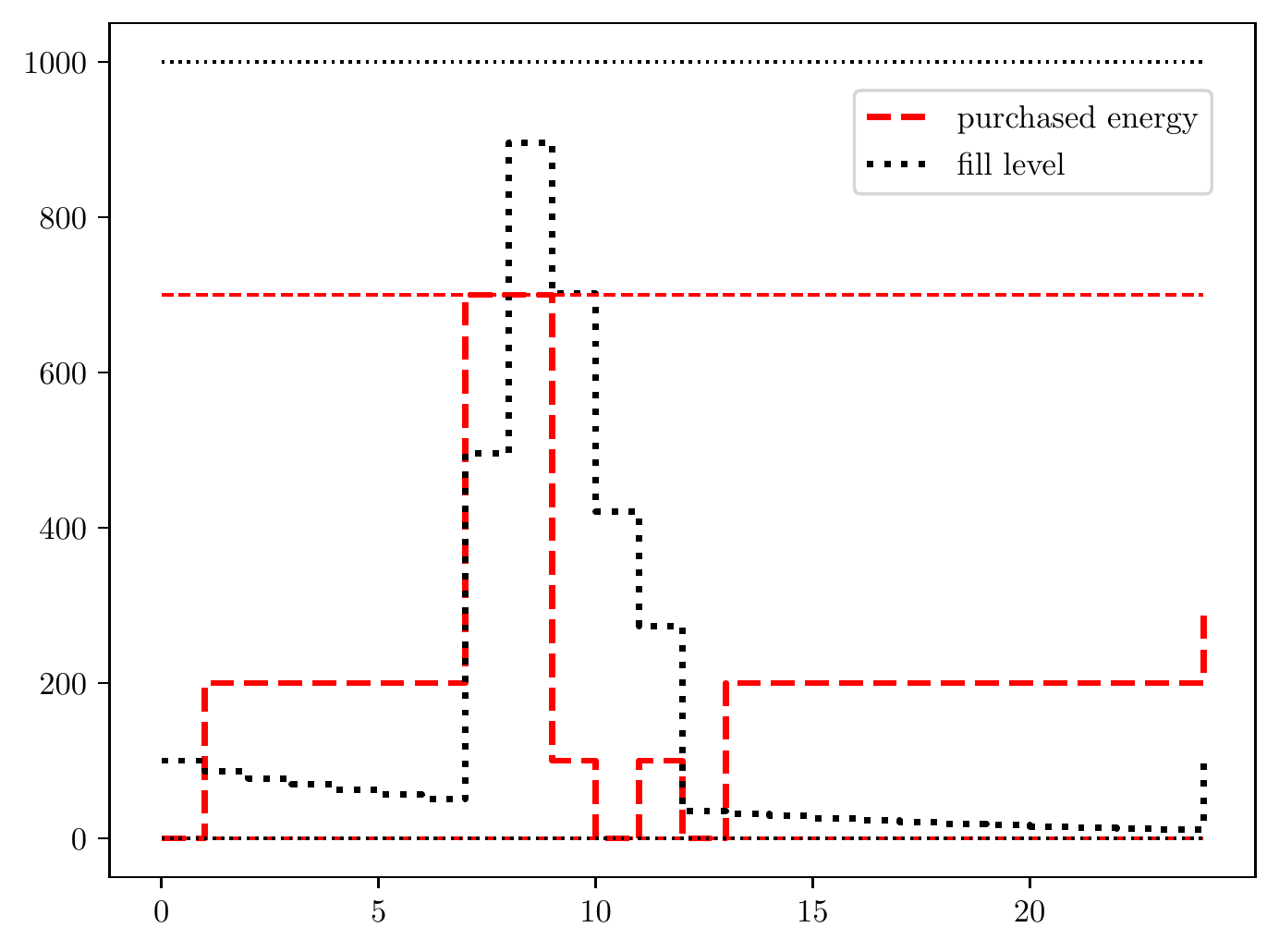}}
    \subfloat[(DP) obj.~value: 1130\,\euro\ ]{\includegraphics[width=0.32\textwidth]{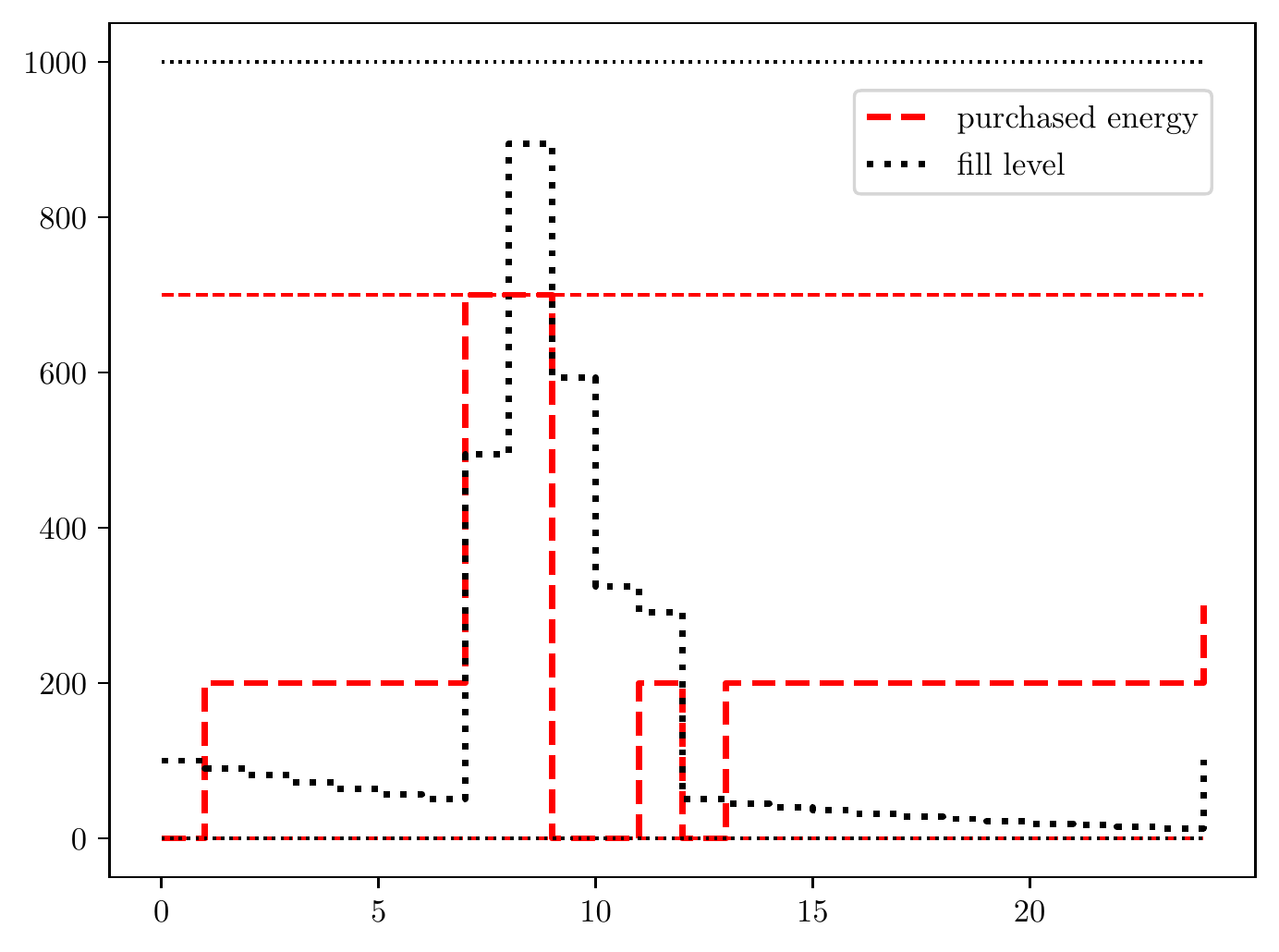}}
    \caption{Comparison of the different optimization methods for 01/07/2018 \new{with a constant energy consumption of $200$\,kWh. The horizontal red and black lines represent the bounds on the capacity and the purchased energy, respectively}.}
    \label{fig:comparison-methods}
\end{figure*}
We implemented the proposed rounding based dynamic programming algorithm, which can be easily adjusted to different use cases. In a simple setup, we compare our method against both linear and mixed integer linear programming, w.r.t.\ the obtained objective function value $z^*$, of the respective approach. We further add some experiments that demonstrate the run-time differences of our method compared to integer programming, as well as a trade-off analysis that contrasts the computed energy costs with the storage capacity. The DP is implemented in Python 3.6, the LP in MATLAB and solved with Gurobi, and the MILP is implemented in Julia and solved with Cbc. All experiments were performed on an Intel(R) Celeron(R) N4000 CPU with 1.10GHz and 8GB main memory.

In our framework there are several parameters to be adjusted. All experiments were performed with a linear energy loss function, however, different types of functions can be applied. Further, we assume a constant energy consumption over the overall time period to obtain interpretable results reflecting the energy costs. According to the nature of the day-ahead market, we allow to purchase energy in steps of $100$\,kWh, whereas the storage is discretized with a step size of $1$\,kWh, \emph{i.e.,} $V_t \in \mathbb{Z}$ for all $t\in T$. A finer or coarser discretization increases the run-time or the rounding error, respectively. In all our numerical experiments we use a time discretization of 1\,h. The storage capacity and the amount of purchasable energy are lower bounded by 0. The maximal fill level is varied in the following experiments, and we restrict the energy that is stored at one time step to be maximum half of the capacity. The efficiency factors are $\eta_\mathrm{in} = 0.9$ and $\eta_\mathrm{out} = 0.95$, the energy loss factor is $\beta=0.1$.

\begin{figure}[t]
    \centering
    \captionsetup{type=figure}
    \includegraphics[width=0.9\textwidth]{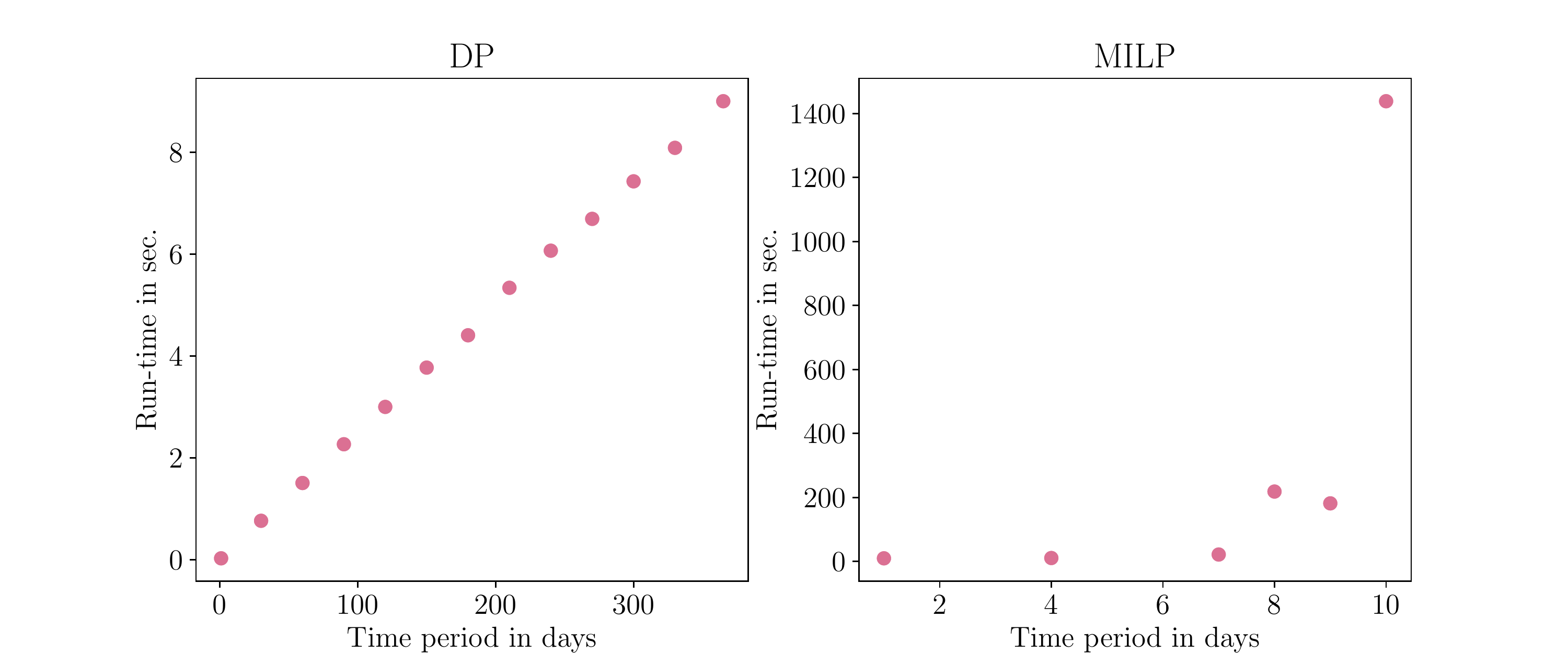}
    \caption{Comparison of the run-time for the DP and the MILP. The run-time is measured in seconds for different time periods. We observe a linearly and exponentially increasing run-time for the DP and MILP, respectively. 
    }
    \label{fig:runtime}

    \subfloat[\new{August 2018}]{\includegraphics[width=0.475\textwidth]{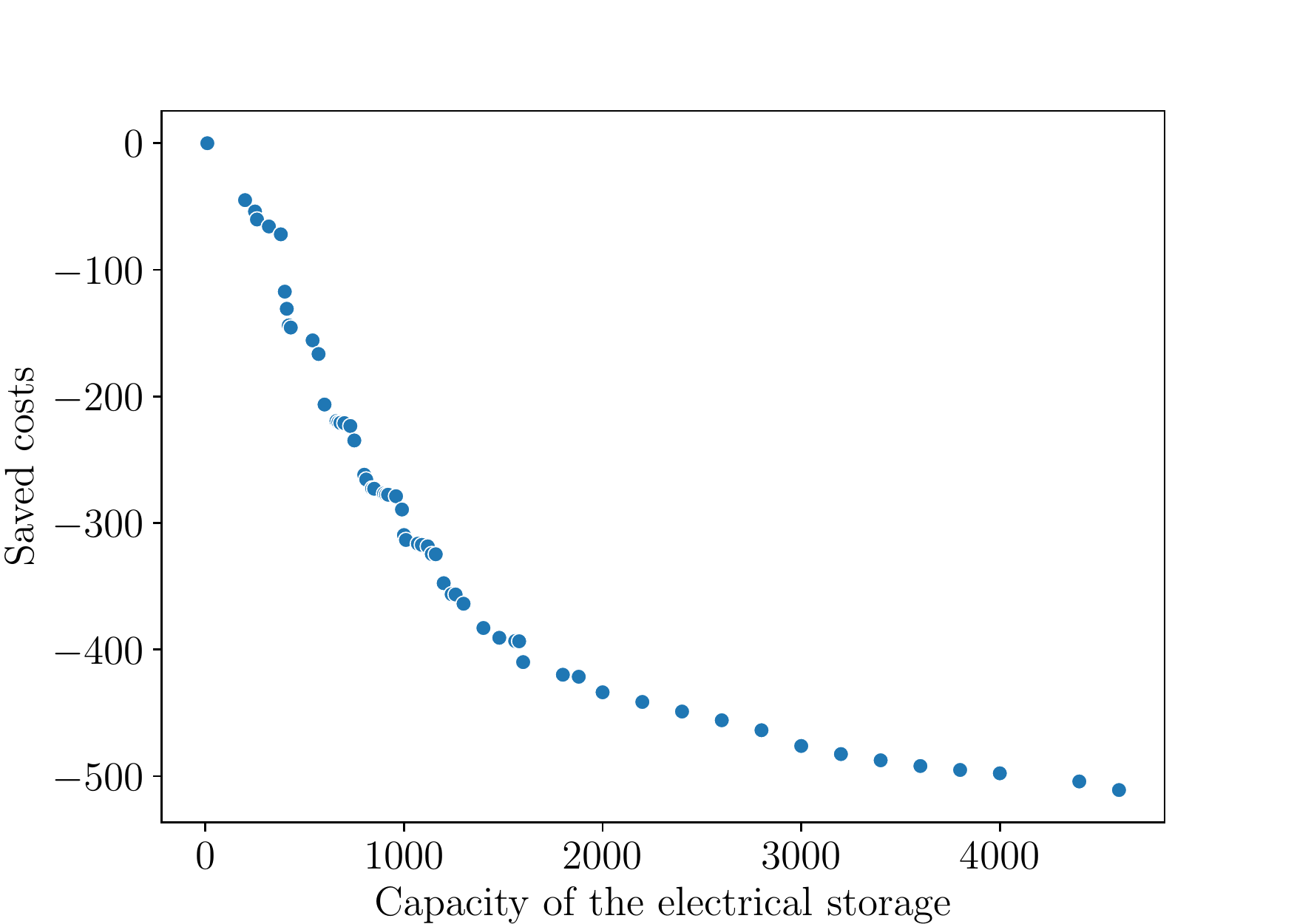}} \hfill
    \subfloat[\new{2018}]{\includegraphics[width=0.475\textwidth]{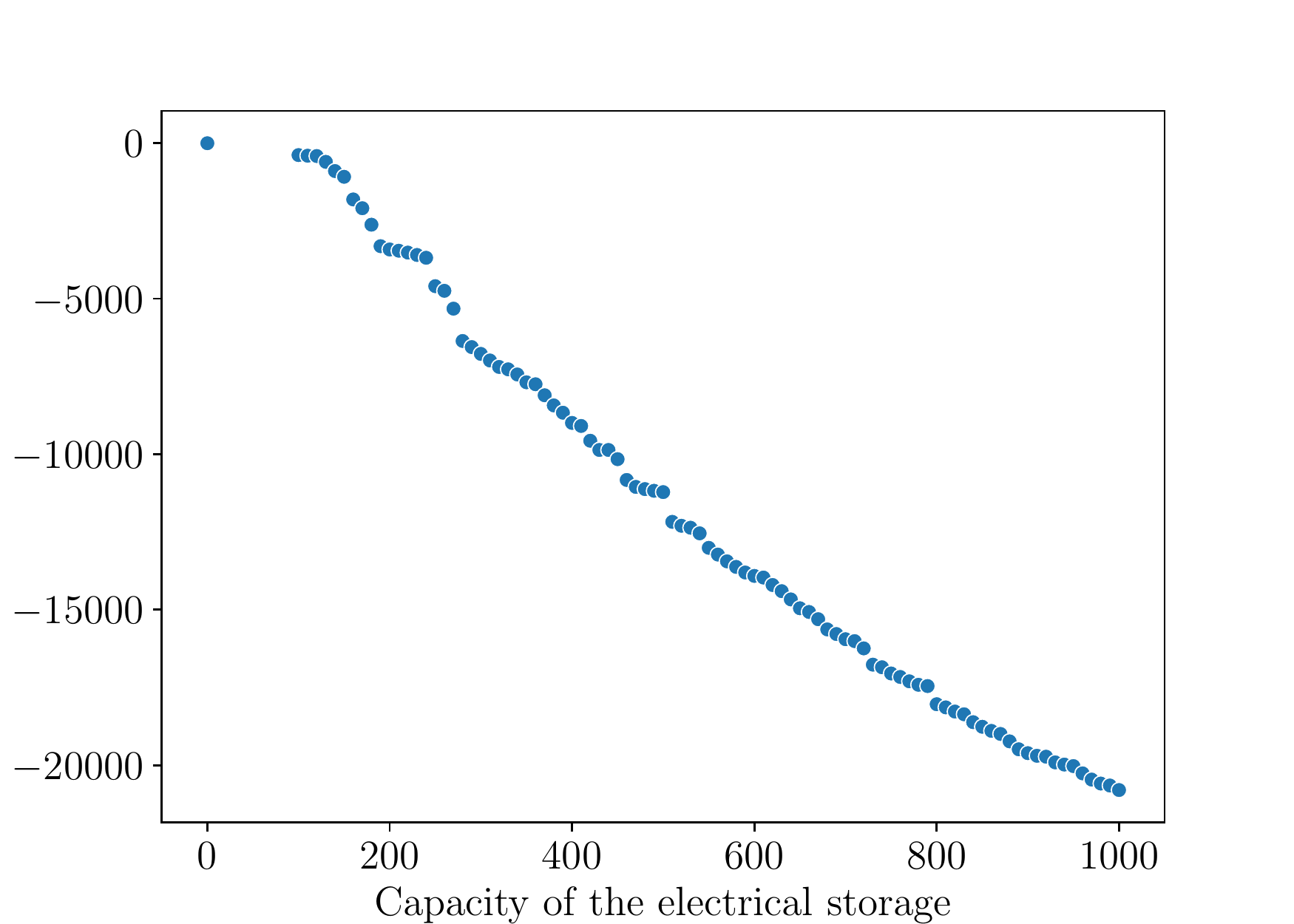}} 
    \caption{Trade-off analysis of energy costs over \new{(a) one month and (b) one year} (in \euro) vs.\ storage capacity (in kWh). The cost is given in comparison to the case without storage\new{, \emph{i.e.}, showing the cost savings}.}
    \label{fig:tradeoff}
\end{figure}

In Table \ref{tab:comparison} we compare the results of our approach against linear and mixed-integer linear programming solutions. Linear programming (LP) achieves the best solutions, since it considers a relaxation of the discrete problem. However, these solutions are not feasible energy inputs from the EPEX Spot market. Our dynamic programming (DP) approach yields only slightly worse solutions compared to the exact optimal solutions of the mixed integer linear programming (MILP) problem (obtained with coin-or/Cbc \cite{forrest18coin}). All solutions are computed with initial and final fill level equal \new{to $100$\,kWh}. A visual example for the method comparison is given in Figure \ref{fig:comparison-methods}. We observe that all three approaches provide qualitatively similar results.

While the run-time of the rounding based dynamic programming algorithm depends not only on the considered time period, but also on the storage capacity, its discretization level and the purchasable energy, the run-time of the MILP is only little impaired by variations of these parameters. However, our approach has the clear advantage compared to the MILP model that the optimization over longer time periods (months/years) or with  finer time discretizations (15\,min/1\,min) is possible. For fixed bounds regarding the storage size and the purchasable energy, its run-time grows only linearly for an increasing number of time steps, while the run-time of the MILP grows exponentially, see Figure \ref{fig:runtime}.


\begin{figure*}[t]
    \centering
    \captionsetup[subfloat]{labelformat=empty}
    \subfloat[]{{\includegraphics[width=0.95\textwidth]{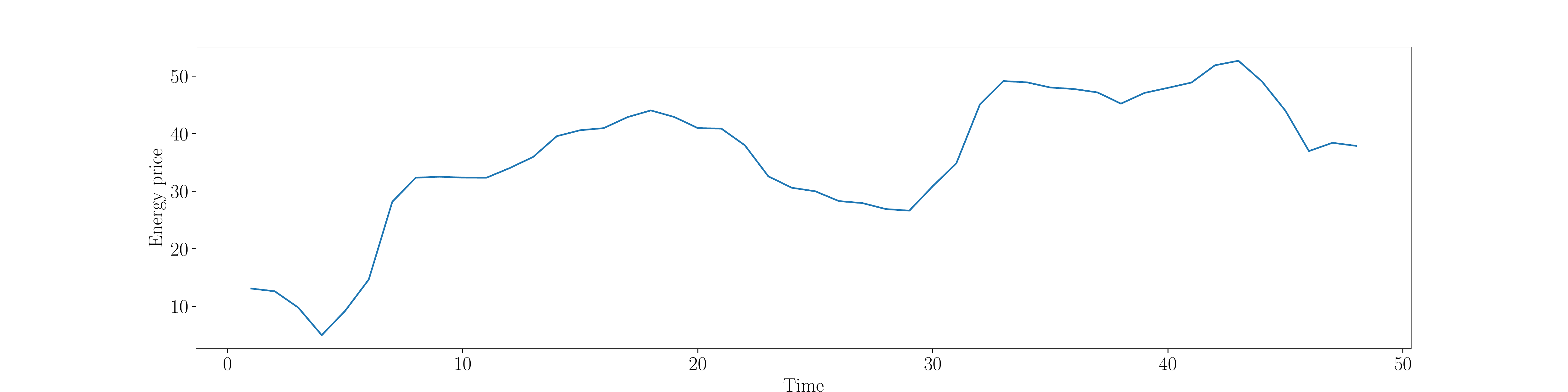}}}\\[-1.5em]
    \captionsetup{type=figure}
    \subfloat[01/08/2018]{\includegraphics[width=0.32\textwidth]{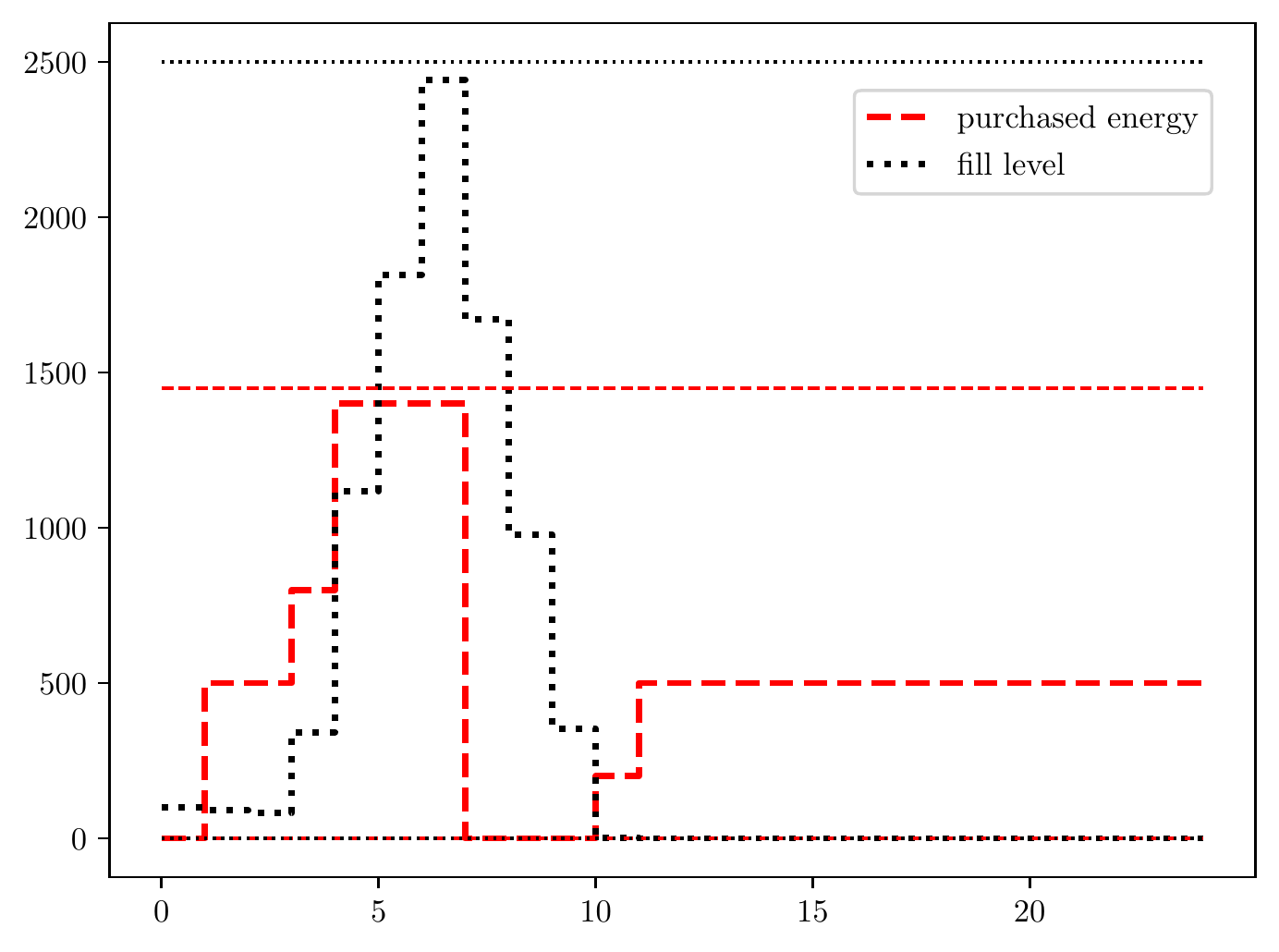}}~
    \subfloat[01/09/2018]{\includegraphics[width=0.32\textwidth]{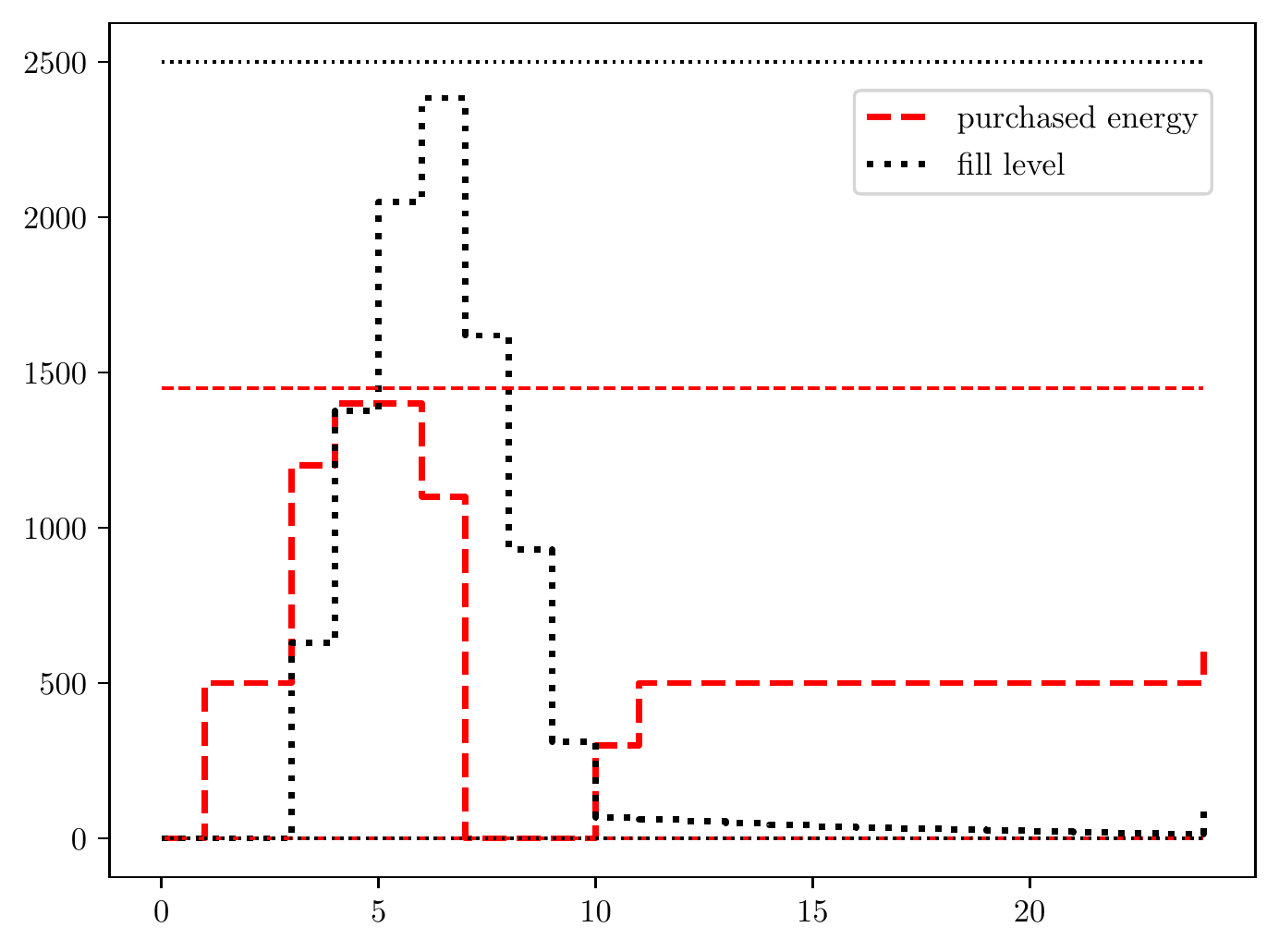}}
    \subfloat[01/08/2018--01/09/2018]{\includegraphics[width=0.32\textwidth]{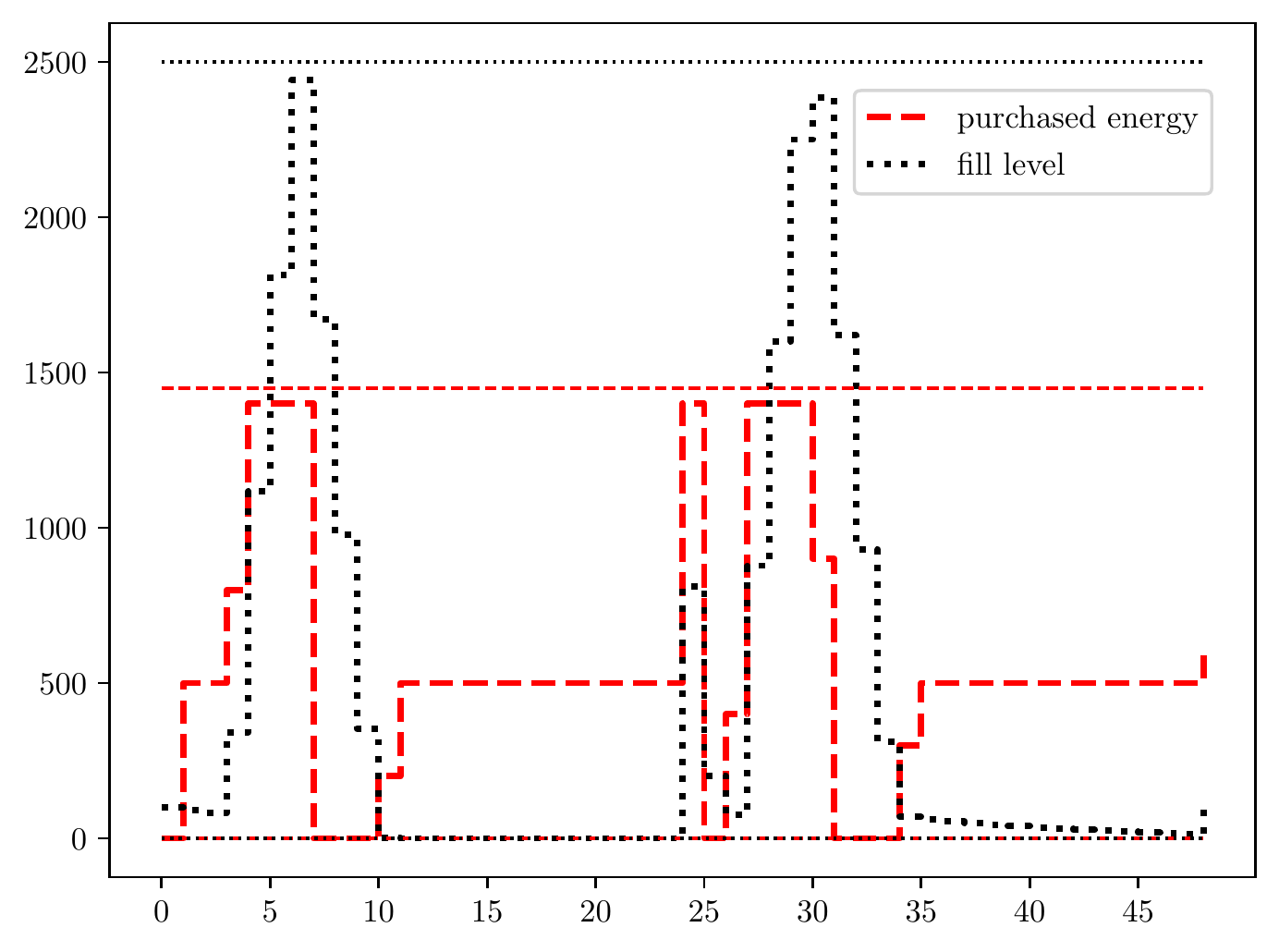}}
    \caption{Optimizing a number of days at once rather than in sequence improves the result. In this example, the costs for optimizing both days separately sum up to 6796\,\euro\ , jointly to 6752\,\euro\,, \new{assuming a constant energy consumption of $200$\,kWh. The horizontal red and black lines represent the bounds on the capacity and the purchased energy, respectively}.}
    \label{fig:comparison-multiple-days}
\end{figure*}

We provide a trade-off analysis (Figure \ref{fig:tradeoff}) where we consider the solutions of our algorithm applied to \new{time frames of one month and one year for varying storage capacities going from $0$\,kWh to $5000$\,kWh and $1000$\,kWh, respectively, in steps of $10$\,kWh.} Based on historical data, previous months or even years can be optimized for several storages and their respective investment costs can be viewed relative to the appropriate energy cost savings. This multiobjective perspective allows to investigate the potential of further investments in storage devices, since it shows the gradual cost reduction induced by the increasing capacity of the storage device. In the case that future energy prices are known either exactly or through forecasting models optimizing several days or weeks jointly improves the result (Figure \ref{fig:comparison-multiple-days}).

\begin{figure}[t]
    \centering
    \captionsetup{type=figure}
    \includegraphics[width=0.9\textwidth]{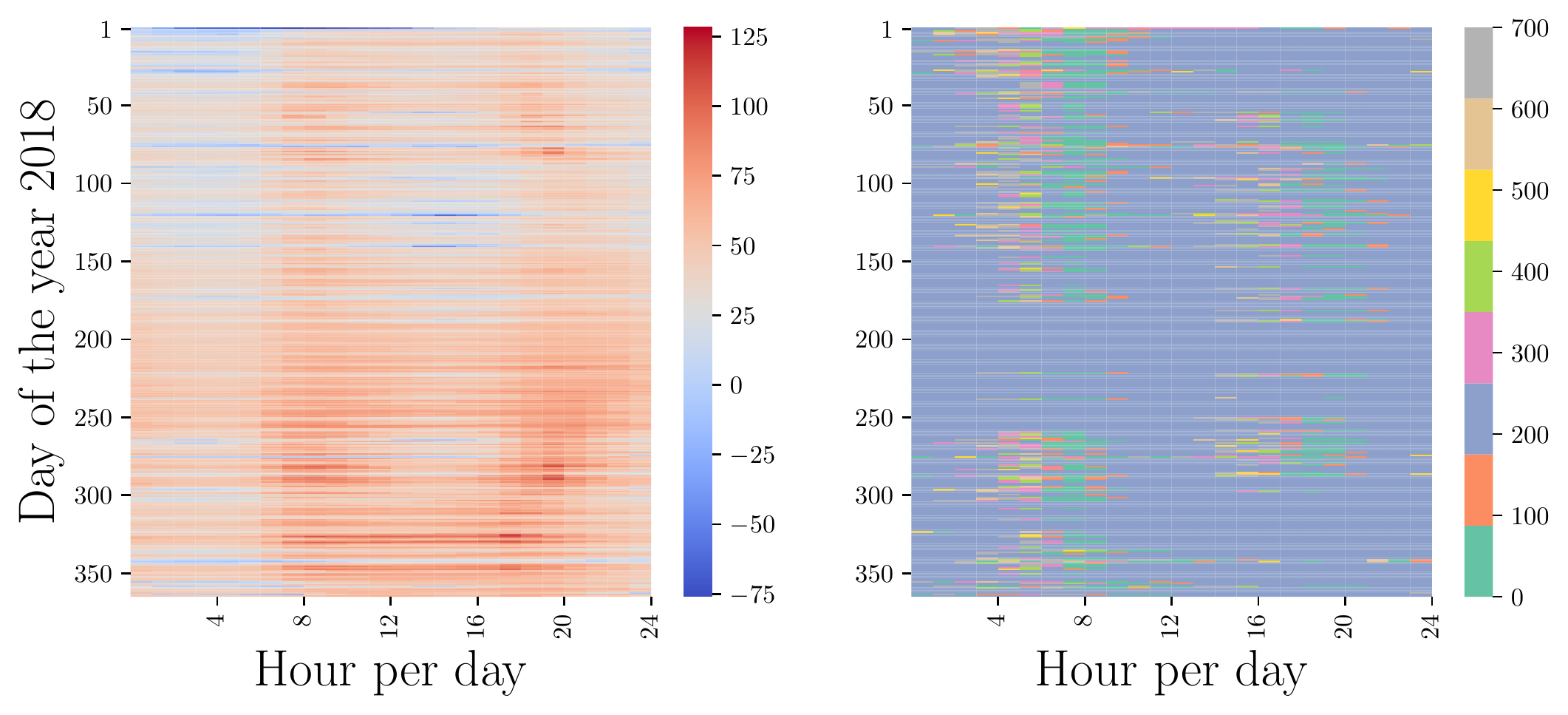}
    \caption{Left: Heatmap that illustrates the hourly energy prices per 100\,kWh for each day of the year 2018. Right: Solution of the DP applied to the year 2018, \emph{i.e.}, purchased energy in kWh for each hour of the year.}
    \label{fig:2018}
\end{figure}

In Figure \ref{fig:2018} we illustrate the solution obtained by the rounding based DP retrospectively optimizing the energy costs over one year, compared to the hourly energy prices. Here, we consider a storage with a maximal fill level of 1000\,kWh and a constant energy consumption of 200\,kWh per hour. We observe that most energy is purchased in the hours before the two energy price peaks. This demonstrates that our algorithm employs the energy storage in order to bridge expensive time periods. 

\section{Conclusions and Outlook}\label{sec:concl}
In this paper, we show that rounding based dynamic programming is an efficient optimization approach for the optimal control of energy storage devices in the presence of volatile costs. In comparison to mixed-integer (linear) programming models the run-time of RBDP is linear in the number of time-steps, which allows us to optimize over larger time periods. The solution of RBDP is thereby a good
approximation of the global optimum obtained by the solution of the MILP model.

The presented computational experiments use simplified load curves. However, it is possible to integrate more complicated load curves and feeding plants, as well as supply from own renewable energy sources. This could be used, for example, to optimize the energy trading of a medium-sized company with its own photovoltaic system and battery storage. In addition, this could also be used to calculate the optimal dimensioning of an energy storage system before the investment.


\begin{backmatter}

\section*{Acknowledgements}
We acknowledge support from the Open Access Publication Fund of the University of Wuppertal. Further we thank Benedikt Dahlmann for the supervision of the master's thesis, on which this work is build upon.

\section*{Funding}
Open Access Publication Fund of the University of Wuppertal.

\section*{Abbreviations}
\begin{tabular}{ll}
DP & dynamic programming\\
LP & linear programming\\
MILP & mixed-integer linear programming\\
MIP & mixed-integer programming\\
RBDP & rounding based dynamic programming
\end{tabular}

\section*{Availability of data and materials}
The prices of the day ahead market are publicly available. The used data are also available together with the code on github: \url{https://github.com/SUhlemeyer/RBDP}.


\section*{Competing interests}
The authors declare that they have no competing interests.


\section*{Authors' contributions}
MS and SU wrote the main part of the paper with input from all authors. SU developed and implemented the rounding based dynamic programming algorithm under supervision of MS and performed the experiments. MS implemented the integer programming problem. BU, supervised by MZ, contributed his knowledge of the energy market and helped with writing the manuscript, especially the introduction and conclusion. All authors contributed to and approved the final manuscript.



\bibliographystyle{vancouver} 
\bibliography{dynamic.bib}      

\end{backmatter}
\end{document}